# Augmented GARCH sequences: Dependence structure and asymptotics

SIEGFRIED HÖRMANN[1]

[1]*Department of Mathematics, University of Utah, 155 S 1400 E RM 233, Salt Lake City, 84112-0090, USA. E-mail: hormann@math.utah.edu*

The augmented GARCH model is a unification of numerous extensions of the popular and widely used ARCH process. It was introduced by Duan and besides ordinary (linear) GARCH processes, it contains exponential GARCH, power GARCH, threshold GARCH, asymmetric GARCH, etc. In this paper, we study the probabilistic structure of augmented GARCH$(1,1)$ sequences and the asymptotic distribution of various functionals of the process occurring in problems of statistical inference. Instead of using the Markov structure of the model and implied mixing properties, we utilize independence properties of perturbed GARCH sequences to directly reduce their asymptotic behavior to the case of independent random variables. This method applies for a very large class of functionals and eliminates the fairly restrictive moment and smoothness conditions assumed in the earlier theory. In particular, we derive functional CLTs for powers of the augmented GARCH variables, derive the error rate in the CLT and obtain asymptotic results for their empirical processes under nearly optimal conditions.

## 1. Introduction

The seminal work of Engle [17] gave a new impact to the theory of time series analysis. Engle introduced the ARCH (autoregressive conditionally heteroscedastic) process, which allows the conditional variance of the time series to change as a function of past observations. In past decades, this model has been widely used in the econometrics literature to describe financial data showing time-varying volatility. Starting with the GARCH model of Bollerslev [8], a great variety of generalizations and extensions of the ARCH model have been introduced and studied. Duan [16] defined a very general model, the so-called augmented GARCH process, which contains most of the known GARCH models as special cases.

For some random variable $Z$, we define $D(Z)$ as the smallest (finite or infinite) interval such that $Z \in D(Z)$ a.s.

**Definition 1.** *Let $g(x)$ and $c(x)$ be real-valued and measurable functions and assume that*

$$\{\varepsilon_k, k \in \mathbb{Z}\} \text{ is an i.i.d. sequence.} \tag{1}$$

*Assume, further, that the stochastic recurrence equation*

$$X_k = c(\varepsilon_{k-1})X_{k-1} + g(\varepsilon_{k-1}), \qquad k \in \mathbb{Z}, \tag{2}$$

*has a strictly stationary solution and assume that*

$$\Lambda: \mathbb{R}^+ \to D(X_0) \text{ is an invertible function.} \tag{3}$$

*Then, an augmented* GARCH$(1,1)$ *process* $\{y_k, -\infty < k < \infty\}$ *is defined by the equation*

$$y_k = \sigma_k \varepsilon_k, \tag{4}$$

*with* $\sigma_k^2 = \Lambda^{-1}(X_k)$.

Necessary and sufficient conditions for the existence of a strictly stationary solution of (2) are given in Aue *et al.* [1] and Duan [16]. In fact, (2) defines a stochastic recurrence equation of the type treated in [9] and [10]. For the convenience of the reader, we state some important results in Section 2 below.

Statistical inference for augmented GARCH sequences requires the study of the asymptotics of the process. For example, detecting structural breaks in GARCH models via CUSUM or MOSUM statistics or determining the limit distribution of the Dickey–Fuller statistics used in unit root testing in AR-GARCH processes requires functional CLTs for partial sums of various functionals of the process. A possible method to prove such results is to use the Markov structure of the model to obtain mixing properties and then to employ the theory of mixing variables to deduce the desired asymptotic results. A detailed study of $\beta$-mixing properties was undertaken by Carrasco and Chen [11]. Mixing properties for ordinary GARCH$(p, q)$ sequences were derived by Basrak *et al.* [3]. An overview of different dependence structures occurring in econometric time series models can be found in Nzé and Doukhan [24]. Besides mixing, they consider, for example, NED (near-epoch dependence), association and $(\psi, \mathcal{L}, \theta)$-weak dependence. An elaborate treatment of the latter concept with applications is presented in Dedecker *et al.* [13]. The great advantage of this approach is that once a good dependence condition is obtained, many asymptotic results (such as functional CLTs) follow directly from a general theory, without any additional work. On the other hand, the class of "ready to use" limit theorems for dependent sequences is limited and the verification of many well-known dependence measures (especially mixing properties) is usually difficult. Hence, to prove the availability of a specific dependence structure often requires fairly restrictive moment and smoothness conditions on the underlying process. For example, in the case of the AR-GARCH model, all existing results on the Dickey–Fuller test, as well as the underlying functional CLT's, require the existence of four moments of the process. Technically, the existence of moments of $y_0$ is a restriction on the functions $c(\cdot)$ and $g(\cdot)$ in (2). Especially in higher order GARCH models, it is very difficult to connect the moments of $y_0$ with the specific model (cf. Ling and McAleer [22]). In this paper, we derive diverse independence properties of GARCH and related models and use these results to directly reduce their asymptotic behavior to that of independent random variables, extending their scope of application and eliminating the restrictive moment and smoothness conditions required by existing theory. This method, adapted from Berkes and Horváth [5], depends on the observation that truncating the infinite series in explicit representations for a GARCH process $\{y_n, n \geq 1\}$ leads to a perturbed process $\{\tilde{y}_n, n \geq 1\}$, whose finite segments $\{\tilde{y}_n, 1 \leq n \leq N\}$ are $m$-dependent, with $m = m(N)$. This fact leads, via simple blocking arguments, to a large class of asymptotic results under optimal conditions. In particular, we will derive functional CLT's for powers of augmented GARCH variables, requiring only the existence of finite variances (Section 4.1). Further, we prove a strong approximation of the empirical process by a two-parameter Gaussian process under logarithmic moment assumptions (Section 4.3). Our method will also yield essentially sharp convergence rates to the normal law (Section 4.2). It is very likely that this approach can be used in many other situations as well.

The structure of our paper is as follows. In Section 2, we formulate preliminary results concerning the existence of a stationary solution to (2) and conditions for $\mathbb{E}|y_0|^p < \infty$ ($p > 0$). In Section 3, we study the dependence structure of augmented GARCH sequences. We deduce an approximation of the original random variables by an $m$-dependent sequence and give estimates for the resulting error. These structural results are applied in Section 4 to derive diverse limit theorems under very mild, or even optimal, assumptions. Finally, Section 5 contains proofs.

## 2. Preliminaries

We give conditions that ensure the existence of a strictly stationary solution of the stochastic recurrence equation (2). Since we assume (1) and (3), this will obviously provide a strictly stationary solution of (4). The following existence theorem is due to Brandt [10]. Let $\log^+ x = \log(x \vee 1)$.

**Theorem 1.** *Assume that* (1) *holds and that* $\mathbb{E}\log^+|g(\varepsilon_0)|$ *and* $\mathbb{E}\log^+|c(\varepsilon_0)|$ *are finite. If*

$$\mathbb{E}\log|c(\varepsilon_0)| < 0, \tag{5}$$

*then for every* $k \in \mathbb{Z}$, *the series*

$$X_k = \sum_{i=1}^{\infty} g(\varepsilon_{k-i}) \prod_{1 \leq j < i} c(\varepsilon_{k-j}) \tag{6}$$

*is convergent with probability one and* $\{X_k, k \in \mathbb{Z}\}$ *is the unique strictly stationary solution of* (2).

If $\{X_k, k \in \mathbb{Z}\}$ is a strictly stationary solution of (2), then it is called *non-anticipative* if $X_k$ is independent of $\sigma(\varepsilon_j, j \geq k)$. Following Bougerol and Picard [9], we call the model (2) *irreducible* if there is no trivial solution $X_k = x \in \mathbb{R}$ such that $x = c(\varepsilon_0)x + g(\varepsilon_0)$ a.s., hence, we exclude the case where $g(\varepsilon_0) = x(1 - c(\varepsilon_0))$ for some $x \in \mathbb{R}$. The following converse to Theorem 1 is a special case of Bougerol and Picard [9], Theorem 2.4.

**Theorem 2.** *Assume that* (1) *holds and that the model* (2) *is irreducible and has a strictly stationary non-anticipative solution. Then, the series in* (6) *converges with probability one and is the unique strictly stationary solution of* (2).

The series representation (6) is essential for our approach since it will provide a construction method for approximating $m$-dependent random variables. From Theorem 2, we learn that whenever a non-anticipative and strictly stationary solution of (2) exists, it has the form (6). Furthermore, Theorem 1 provides easy conditions for the convergence of this series. Hence, instead of requiring the existence of a solution of (2), we will hereafter assume a priori that the series in (6) is convergent and that

$$\Lambda(\sigma_k^2) = \sum_{i=1}^{\infty} g(\varepsilon_{k-i}) \prod_{1 \leq j < i} c(\varepsilon_{k-j}). \qquad (7)$$

We shall abbreviate this compound assumption as simply "(7) holds". Also, note that (7) is always a formal solution of (2).

The next theorem is only a slight modification of a result of Aue *et al.* [1].

**Theorem 3.** *Assume that* (1) *and* (7) *hold. If, for some $\mu > 0$,*

$$\mathbb{E}|g(\varepsilon_0)|^\mu < \infty \quad and \quad \mathbb{E}|c(\varepsilon_0)|^\mu < 1, \qquad (8)$$

*then $\mathbb{E}|\Lambda(\sigma_0^2)|^\mu < \infty$. Conversely, if*

$$c(\varepsilon_0) \geq 0 \quad and \quad g(\varepsilon_0) \geq 0 \qquad (9)$$

*hold, then* (8) *is necessary in order that $\mathbb{E}|\Lambda(\sigma_0^2)|^\mu < \infty$.*

**Remark 1.** The sufficiency of (8) for $\mathbb{E}|\Lambda(\sigma_0^2)|^\mu < \infty$ is also given in Carrasco and Chen [11], but under additional assumptions. They require $\mu \in \{1, 2, \ldots\}$, $\mathbb{E}\varepsilon_0 = 0$, $\mathbb{E}\varepsilon_0^2 = 1$ and a continuous density $p(\cdot)$ of $\varepsilon_0$ such that $p(x) > 0$ for all $x \in \mathbb{R}$.

In [11] and [16], several examples of augmented GARCH(1,1) processes appearing in the literature are given. We refer to these articles for an overview. It should be noted that most commonly known examples (as well as all examples given in [11, 16]) treat so-called *polynomial* or *exponential* GARCH processes. We say that an augmented GARCH sequence is a polynomial (resp., exponential) GARCH process if $\Lambda(x) = x^\delta$ with $\delta > 0$ (resp., $\Lambda(x) = \log x$). This special form of $\Lambda$ is motivated from a Box–Cox transformation of the observation. In order to avoid cumbersome calculations for only a small gain, in this paper, we shall solely deal with these two types. Our approach will work with more general $\Lambda$ if we require certain smoothness conditions as, for example, in Aue *et al.* [1]. However, in a generalized formulation, the price will be a loss of accuracy of our results. The main advantage of a special choice of $\Lambda$ is that we can very precisely transfer different properties which are derived for $\Lambda(\sigma_k^2)$ to the volatility process $\{\sigma_k^2\}$.

The following corollary to Theorem 3 characterizes the existence of $\mathbb{E}|y_0|^{2\nu}$ ($\nu > 0$) for some polynomial GARCH sequence. (See also Ling and McAleer [23].)

**Corollary 1.** *Assume that* (1), (4), (7), (9) *and* $\mathbb{E}|\varepsilon_0|^{2\nu} < \infty$ *hold. Assume, further, that* $\Lambda(x) = x^\delta$ ($\delta > 0$). *Then, for some* $\nu > 0$, *we have* $\mathbb{E}|y_0|^{2\nu} < \infty$ *if and only if*

$$\mathbb{E}|c(\varepsilon_0)|^{\nu/\delta} < 1 \quad and \quad \mathbb{E}|g(\varepsilon_0)|^{\nu/\delta} < \infty. \tag{10}$$

For polynomial GARCH processes, the left-hand side of (2) is always positive and thus it is natural to assume (9). Note, however, that $c(\varepsilon_0) \geq 0$ and $g(\varepsilon_0) \geq 0$ are not required for $\Lambda(\sigma_k^2) \geq 0$. A non-trivial example is where $g(x) = 1$ and the distribution of $c(\varepsilon_0)$ is concentrated on the interval $[-1/2, 0]$. It follows from the definition of $\Lambda$ that $\mathbb{E}|y_0|^{2\nu} < \infty$ if and only if $\mathbb{E}|\varepsilon_0|^{2\nu} < \infty$ and $\mathbb{E}|\Lambda(\sigma_0^2)|^{\nu/\delta} < \infty$, which, in view of Theorem 3, proves the corollary.

Proposition 1 below gives a moment criterion for exponential GARCH processes. For the special case $c(\varepsilon_0) = c$ with $|c| < 1$, we refer to He *et al.* [21].

**Proposition 1.** *Assume that* (1), (4), (7) *and* $\mathbb{E}|\varepsilon_0|^{2\mu} < \infty$ *hold. Assume, further, that* $\Lambda(x) = \log x$. *If*

$$|c(\varepsilon_0)| \leq c < 1 \quad and \quad \mathbb{E}e^{\mu|g(\varepsilon_0)|} < \infty, \tag{11}$$

*then* $\mathbb{E}|y_0|^{2\mu} < \infty$. *On the other hand, if* (9) *holds, then*

$$P(|c(\varepsilon_0)| \leq 1) = 1 \tag{12}$$

*and the second assumption in* (11) *are necessary to assure* $\mathbb{E}|y_0|^{2\mu} < \infty$.

The sufficiency part is implicit in [1]. To prove the other direction, we observe that by (7) and the non-negativity of $c(\varepsilon_0)$ and $g(\varepsilon_0)$, it follows that

$$\mathbb{E}|y_0|^{2\mu} = \mathbb{E}|\varepsilon_0|^{2\mu}\mathbb{E}\exp(\mu\Lambda(\sigma_0^2)) \geq \mathbb{E}|\varepsilon_0|^{2\mu}\mathbb{E}\exp(\mu g(\varepsilon_0)) \tag{13}$$

and hence $\mathbb{E}e^{\mu g(\varepsilon_0)} < \infty$ is necessary. Further, it is clear from (13) that $\Lambda(\sigma_0^2)$ must have moments of all orders. If we assume that (12) does not hold, then there is some $\delta > 0$ and some $\alpha > 0$ such that $P\{c(\varepsilon_0) \geq 1 + \delta\} = \alpha$. Consequently, $\mathbb{E}c(\varepsilon_0)^p \geq \alpha(1+\delta)^p > 1$ for sufficiently large $p$. By (7) and the assumption $c(\varepsilon_0), g(\varepsilon_0) \geq 0$, we have

$$\mathbb{E}\Lambda(\sigma_0^2)^p \geq \sum_{i=1}^\infty \mathbb{E}g(\varepsilon_0)^p(\mathbb{E}c(\varepsilon_0)^p)^{i-1} = \infty,$$

which contradicts the fact that $\mathbb{E}\Lambda(\sigma_0^2)^p < \infty$ for any $p > 0$.

## 3. Dependence structure

In this section, we show that augmented GARCH sequences $(y_k)$ can be closely approximated by $m$-dependent random variables $(y_{km})$. The error of the approximation will be measured by the $L_2$-distance

$$\|y_k - y_{km}\|_2 = (\mathbb{E}|y_k - y_{km}|^2)^{1/2}$$

and by probability inequalities for

$$P(|y_k - y_{km}| > \epsilon).$$

Under the present assumptions, such types of inequalities are relatively easy to obtain, yet they have very strong consequences. In Section 4, we show that the closeness of the $y_k$'s to independence allows the study of asymptotic properties with classical methods and the derivation of sharp limit theorems via a unified approach. Since, in many applications, it is not the original observations that are analyzed, but the transformed variables $\eta_k = f(y_k)$, where $f$ comprises the functions

$$f(x) = |x|^\nu \quad \text{or} \quad f(x) = \text{sign}(x) \cdot |x|^\nu \qquad (\nu > 0), \tag{14}$$

we shall derive our results for $(\eta_k)$. The basic idea for the approximation is this. If (7) holds, then

$$\sigma_k^2 = \Lambda^{-1}\left(\sum_{1 \leq i < \infty} g(\varepsilon_{k-i}) \prod_{1 \leq j < i} c(\varepsilon_{k-j})\right).$$

We now define

$$y_{km} = \varepsilon_k \sigma_{km}, \tag{15}$$

where $\sigma_{km}^2$ is the solution of

$$\Lambda(\sigma_{km}^2) = \sum_{1 \leq i \leq m} g(\varepsilon_{k-i}) \prod_{1 \leq j < i} c(\varepsilon_{k-j}). \tag{16}$$

It follows that $(y_{km})$ defines an $m$-dependent sequence. Usually, the series in (16) converges very fast and hence considering the finite sums will only cause a small error.

**Lemma 1.** *Assume that* (1), (4), (7) *and* (9) *hold. Assume, further, that* $\Lambda(x) = x^\delta$. *Let $f$ be given as in* (14) *and set*

$$\eta_k = f(y_k) \quad \text{and} \quad \eta_{km} = f(y_{km}). \tag{17}$$

If $\mathbb{E}|\eta_0|^2 < \infty$, then there are constants $C_1 > 0$ and $0 < \varrho < 1$ such that

$$\|\eta_k - \eta_{km}\|_2 \leq C_1 \varrho^m \qquad (m \geq 1).$$

**Lemma 2.** *Assume that* (1), (4), (7) *and* (11) *hold. Assume, further, that* $\Lambda(x) = \log x$. *Let* $\eta_k$ *and* $\eta_{km}$ *be defined as in* (17) *and assume that* $\mathbb{E}|\eta_0|^2 < \infty$. *If* $\mu > \nu$, *then there is some* $C_2 > 0$ *such that*

$$\|\eta_k - \eta_{km}\|_2 \leq C_2 c^m \qquad (m \geq 1). \tag{18}$$

*(Here, $\mu$ and $\nu$ come from* (11) *and* (14) *resp.).*

**Remark 2.** Under the assumption $\mathbb{E}|\eta_0|^p < \infty$ for $p \geq 1$, some trivial changes in the proofs also yield the exponential decrease of the $L_p$-error.

The conclusion of Lemmas 1 and 2 is that under some minor regularity assumptions – providing manageable conditions for the existence of $\mathbb{E}|\eta_0|^2$ (see Corollary 1 and Proposition 1) – the $L_2$ approximation error via the perturbed sequence decays exponentially fast. An immediate consequence is the short-memory behavior of the sequence $(\eta_k)$, which can be characterized by the convergence of

$$\sum_{k=1}^{\infty} |\mathrm{Cov}(\eta_0, \eta_k)| < \infty.$$

In fact, for $m = m(k) = k - 1$, we get, by the independence of $\eta_0$ and $\eta_{km}$, that $|\mathrm{Cov}(\eta_0, \eta_k)| = |\mathrm{Cov}(\eta_0, \eta_k - \eta_{km})| \leq 4\|\eta_0\|_2 \|\eta_k - \eta_{km}\|_2$.

Having established the $L_p$-error of the approximating r.v.'s, one can immediately estimate $P(|\Lambda(\sigma_k^2) - \Lambda(\sigma_{km}^2)| > \epsilon)$ and $P(|y_k - y_{km}| > \epsilon)$ by using the Markov inequality. However, with a little more effort, one can drastically weaken the associated moment requirements. Observe that under (7),

$$|\Lambda(\sigma_k^2) - \Lambda(\sigma_{km}^2)| = \prod_{i=1}^{m} |c(\varepsilon_{k-i})| |\Lambda(\sigma_{k-m}^2)|.$$

If (5) holds, then we can choose $c_1 > 0$ such that $-c_1 - \mathbb{E}\log|c(\varepsilon_0)| > 0$. From the stationarity of $\Lambda(\sigma_k^2)$, we get

$$P\left(\prod_{i=1}^{m} |c(\varepsilon_{k-i})| |\Lambda(\sigma_{k-m}^2)| > e^{-c_1 m/2}\right)$$
$$\leq P\left(\sum_{i=1}^{m} (\log|c(\varepsilon_{k-i})| - \mathbb{E}\log|c(\varepsilon_0)|) > (-c_1 - \mathbb{E}\log|c(\varepsilon_0)|)m\right) + P(|\Lambda(\sigma_0^2)| > e^{c_1 m/2}).$$

Of course, the first term in the previous line only makes sense if $\mathbb{E}\log|c(\varepsilon_0)| > -\infty$. These simple observations are summarized in the next lemma.

**Lemma 3.** *Assume that* (1) *and* (7) *hold. Assume, further, that*

$$-\infty < \mathbb{E}\log|c(\varepsilon_0)| < 0 \tag{19}$$

*and set*

$$T_m = \sum_{i=1}^{m}(\log|c(\varepsilon_i)| - \mathbb{E}\log|c(\varepsilon_i)|).$$

*There are then constants $C_3, C_4 > 0$ which do not depend on $m$ such that for sufficiently small $\alpha > 0$, we have*

$$P(|\Lambda(\sigma_k^2) - \Lambda(\sigma_{km}^2)| > \mathrm{e}^{-\alpha m}) \leq P(T_m > C_3 m) + P(|\Lambda(\sigma_0^2)| > \mathrm{e}^{C_4 m}).$$

Lemma 3 shows that the error rate depends on the tail probabilities of $T_m$ and $|\Lambda(\sigma_0^2)|$. Both can be obtained by classical methods for independent random variables. Clearly, these probabilities are determined by the specific distributions of $c(\varepsilon_0)$ and $g(\varepsilon_0)$. The next two lemmas give more concrete estimates under different moment assumptions.

**Lemma 4.** *Assume that the conditions of Lemma 3 hold. If for some $\mu > 2$,*

$$\mathbb{E}|\log|c(\varepsilon_0)||^\mu < \infty \quad \text{and} \quad \mathbb{E}(\log^+|g(\varepsilon_0)|)^\mu < \infty, \tag{20}$$

*then for a sufficiently small $\alpha > 0$, there is some constant $C_5$ such that for every $m \geq 1$,*

$$P(|\Lambda(\sigma_k^2) - \Lambda(\sigma_{km}^2)| > \mathrm{e}^{-\alpha m}) \leq C_5 m^{(2-\mu)/2}.$$

**Lemma 5.** *Assume that the conditions of Lemma 3 hold. If for some $\mu > 0$,*

$$\mathbb{E}|c(\varepsilon_0)|^\mu < 1 \quad \text{and} \quad \mathbb{E}|g(\varepsilon_0)|^\mu < \infty, \tag{21}$$

*then for a sufficiently small $\alpha > 0$, there are constants $C_6, \rho > 0$ such that for every $m \geq 1$,*

$$P(|\Lambda(\sigma_k^2) - \Lambda(\sigma_{km}^2)| > \mathrm{e}^{-\alpha m}) \leq C_6 \mathrm{e}^{-\rho m}.$$

Finally, we consider the estimation of $P(|y_k - y_{km}| > \epsilon)$. As in Lemmas 1–2, we will be a little bit more general, considering the transformed variables $\eta_k$ and distinguisting between polynomial and exponential GARCH.

**Lemma 6.** *Assume that* (1), (4), (7), (9) *and* (19) *hold and that* (20) *holds with some $\mu > 2$. Assume, further, that $\Lambda(x) = x^\delta$ and let $\eta_k$ and $\eta_{km}$ be given as in* (17). *Then, for sufficiently small $\alpha > 0$, there is some constant $C_7$ such that*

$$P(|\eta_k - \eta_{km}| > \mathrm{e}^{-\alpha m}) \leq P(|\varepsilon_0| \geq \mathrm{e}^{\alpha m}) + C_7 m^{(2-\mu)/2} \qquad (m \geq 1).$$

**Lemma 7.** *Assume that* (1), (4), (7) *hold and that* (21) *holds with some* $\mu > 0$. *Assume, further, that* $\Lambda(x) = \log x$ *and let* $\eta_k$ *and* $\eta_{km}$ *be given as in* (17). *Then, for sufficiently small* $\alpha > 0$, *there is some constant* $C_8$ *such that*

$$P(|\eta_k - \eta_{km}| > \mathrm{e}^{-\alpha m}) \leq P(|\varepsilon_0| \geq \mathrm{e}^{\alpha m}) + C_8 m^{-\mu} \qquad (m \geq 1).$$

## 4. Applications

### 4.1. The functional central limit theorem

Statistical inference based on GARCH sequences often requires the establishment of functional central limit theorems for partial sums processes like

$$S_n(t) = \frac{1}{\sqrt{n}} \sum_{1 \leq i \leq nt} (f(y_i) - Ef(y_0))$$

or

$$Z_n(t) = \frac{1}{\sqrt{n}} \sum_{1 \leq i \leq nt} (f(\sigma_i) - Ef(\sigma_0)).$$

For example, to derive the asymptotic distribution of the CUSUM or MOSUM statistics applied in the theory of change-point detection, an FCLT is needed. Also, the determination of the asymptotic distribution of the Dickey–Fuller statistic for the unit root test in an AR-GARCH model requires an FCLT. Several authors have obtained functional limit theorems for $S_n(t)$ with $f(x) = x$, $f(x) = |x|$ or $f(x) = x^2$ under various conditions. Most FCLTs in the literature assume at least four moments of the GARCH variables (cf. Davidson [12], Hansen [19]) instead of the more desirable assumption of solely a finite variance. Berkes *et al.* [6] noted that in the special case $f(x) = x$ and $\mathbb{E}\varepsilon_0 = 0$, the envisaged result follows from an FCLT for ergodic martingale difference sequences (cf. Billingsley [7], Theorem 23.1) under the optimal condition $\mathbb{E}y_0^2 < \infty$. However, if $\mathbb{E}y_0 \neq 0$ or if we are, for example, interested in the squared GARCH sequence, then the martingale structure no longer applies. Giraitis *et al.* [18] point out the important association property of the ARCH($\infty$) model, which, when specialized to the GARCH case, yields an FCLT for the sequence $y_k^2$ under $\mathbb{E}y_k^4 < \infty$. Since association is preserved by non-decreasing transformations, an FCLT for the sequences $(|y_k|^\nu)$ follows if one shows that $\sum_{k=-\infty}^{\infty} \operatorname{Cov}(|y_0|^\nu, |y_k|^\nu) < \infty$. Although the convergence of the last series is implied by Lemmas 1 and 2, this approach will not lead to the desired goal here. Due to the general form of the functions $c(x)$ and $g(x)$ in (2), it is not clear whether the association property of the independent $\varepsilon_k$'s is inherited as it is in the ARCH($\infty$) case. Also, for the variables $\operatorname{sign}(y_k) \cdot |y_k|^\nu$, this method seems not to be applicable. A further weak dependence condition has been proposed by Doukhan and Louhichi [14] (see also Doukhan and Wintenberger [15]), the so-called $(\theta, \mathcal{L}, \psi)$-weak dependence. Roughly speaking, a sequence is $(\theta, \mathcal{L}, \psi)$-dependent if for all $n \geq 1$ and all bounded Lipschitz functions $h, k : \mathbb{R}^n \to \mathbb{R}$, the covariances $\mathrm{Cov}(h(\text{'past'}), k(\text{'future'})) \leq \psi(\cdot)\theta(\cdot)$. Here, $\psi$ and $\theta$ are functions, depending on the choice of 'past' and 'future'. If the covariance relation holds with proper $\psi$ and $\theta$, this structure can be used to verify an FCLT (see Nzé and Doukhan [24], Section 5.1.6). Berkes *et al.* [6], Theorem 2.9, applied this method to derive the FCLT for the squares of some GARCH sequence under the condition $\mathbb{E}|y_0|^\kappa < \infty$ for some $\kappa > 8$.

Here, we will obtain the FCLT for the partial sum processes $S_n(t)$ (resp., $Z_n(t)$) by assuming only the necessary assumption $Ef(y_0)^2 < \infty$. Let $\stackrel{\mathcal{D}[0,1]}{\longrightarrow}$ denote weak convergence of a random function in the Skorokhod space $\mathcal{D}[0,1]$.

**Theorem 4.** *Assume that* (1), (4), (7) *and* (9) *hold. Assume, further, that* $\Lambda(x) = x^\delta$ *($\delta > 0$) and that $f$ is given as in* (14). *If* $\mathbb{E}f(y_0)^2 < \infty$, *then*

$$\tau^2 = \mathrm{Var}\, f(y_0) + 2 \sum_{1 \leq k < \infty} \mathrm{Cov}(f(y_0), f(y_k)) \tag{22}$$

*is convergent and*

$$S_n(t) \stackrel{\mathcal{D}[0,1]}{\longrightarrow} \tau W(t),$$

*where* $\{W(t), 0 \leq t \leq 1\}$ *is a Brownian motion. An analogous result is valid for* $Z_n(t)$.

In particular, when the $y_k$ are GARCH$(1,1)$ variables, we obtain FCLTs

$$\frac{1}{\sqrt{n}} \sum_{k=1}^{nt}(y_k - \mathbb{E}y_k) \stackrel{\mathcal{D}[0,1]}{\longrightarrow} \tau_1 W(t)$$

or

$$\frac{1}{\sqrt{n}} \sum_{k=1}^{nt}(\sigma_k^2 - \mathbb{E}\sigma_k^2) \stackrel{\mathcal{D}[0,1]}{\longrightarrow} \tau_2 W(t) \tag{23}$$

(here, $\tau_1^2, \tau_2^2$ are the corresponding variances arising from (22)) under the necessary conditions $\mathbb{E}y_k^2 < \infty$ and $\mathbb{E}\sigma_k^4 < \infty$, respectively.

**Theorem 5.** *Assume that* (1), (4), (7) *and* (11) *hold. Assume, further, that* $\Lambda(x) = \log x$ *and let $f$ be given as in* (14). *If* $\mathbb{E}f(y_0)^2 < \infty$ *and* $\mu > \nu$, *then the proposition of Theorem 4 holds. (Here, $\mu$ and $\nu$ stem from* (11) *and* (14).*)*

**Proof of Theorems 4 and 5.** Set $\xi_k = \eta_k - \mathbb{E}\eta_0$ and $\xi_{km} = \eta_{km} - \mathbb{E}\eta_0$, where $\eta_k$ and $\eta_{km}$ are given as in (17). Obviously, $\xi_k = g(\ldots, \varepsilon_{k-1}, \varepsilon_k)$, where $g$ is some measurable mapping from the space of infinite sequences into $\mathbb{R}$. According to Theorem 21.1 of Billingsley [7], for the proof, it is enough to find measurable mappings $g_m$ from $\mathbb{R}^m$ into $\mathbb{R}$ such that

$$\sum_{1 \leq m < \infty} \|\xi_0 - \zeta_{0m}\|_2 < \infty, \tag{24}$$

where $\zeta_{0m} = g_m(\varepsilon_{-m+1}, \varepsilon_{-m+2}, \ldots, \varepsilon_0)$. By Lemma 1, $\xi_{0m}$ satisfies this requirement and thus setting $\zeta_{0m} = \xi_{0m}$ yields the proof of Theorem 4. The proof of Theorem 5 is the same. □

### 4.2. Rates of convergence in the central limit theorem

The functional CLT developed in Section 4.1 clearly implies the ordinary CLT for the partial sums and raises the question of the normal approximation error. If $f$ is given as in (14), we will obtain the rate of convergence of

$$S_n = f(y_1) + \cdots + f(y_n) - n\mathbb{E}f(y_1)$$

to the normal distribution provided that $m_3 = \mathbb{E}|f(y_1)|^3 < \infty$. For notational convenience, we will write $S_n$ instead of $S_n(f)$ and we set $B_n^2 = B_n(f)^2 = \operatorname{Var} S_n(f)$ and

$$\beta_n^2 = \beta_n^2(f) = \operatorname{Var} f(y_0) + 2 \sum_{j=1}^{n-1} (1 - j/n) \operatorname{Cov}(f(y_1), f(y_{j+1})),$$

that is, $\beta_n^2 = B_n^2/n$.

**Theorem 6.** *Assume that the conditions of Theorem 4 hold with the added assumption $\mathbb{E}|f(y_0)|^3 < \infty$. Then,*

$$\beta := \lim \beta_n \ exists \tag{25}$$

*and if $\beta > 0$, there is some $C > 0$ such that*

$$|P\{S_n < xB_n\} - \Phi(x)| \leq C \frac{(\log n)^2}{\sqrt{n}} \qquad \text{for all } n \geq 2 \text{ and } x \in \mathbb{R}.$$

*The constant $C$ may depend on $f$, $\Lambda$, $c$, $g$ and the law of $\varepsilon_0$.*

**Theorem 7.** *Assume that the conditions of Theorem 5 hold with the added assumptions $\mathbb{E}|f(y_0)|^3 < \infty$ and $\mu > 3\nu/2$. Then, the proposition of Theorem 6 holds.*

The additional assumptions in Theorems 6–7, compared with Theorems 4–5, arise from the requirement $\mathbb{E}|f(y_1)|^3 < \infty$, which is the classical assumption in the context of Berry–Esseen bounds. The existence of the limit in (25) follows from Theorems 4–5. The rate $(\log n)^2 n^{-1/2}$ coincides with that obtained by Tihomirov [27], Theorem 2, for sequences which are $\beta$-mixing with geometric rate. The proofs of Theorems 6–7 are given in Section 5.2. The main ingredients are a Berry–Esseen bound for $m$-dependent sequences (Tihomirov [27], Theorem 5) and Lemmas 1–2. Again, it becomes clear that $m$-dependence, rather than mixing, is the crucial structural property required in order to study the asymptotics of augmented GARCH variables.

### 4.3. Asymptotics of the empirical process

In this section, we study the empirical process of augmented GARCH sequences. For this purpose, we define

$$R(s,t) = \sum_{1 \leq k \leq t} (I\{F(y_k) \leq s\} - s),$$

where $F(x) = P(y_0 \leq x)$. We will derive an almost sure approximation theorem for $R(s,t)$ by a two-parameter Gaussian process $K(s,t)$, assuming only the existence of *logarithmic* moments of the $\varepsilon_k$ and $c(\varepsilon_k)$, $g(\varepsilon_k)$. In the sequel, we shall use $Y_k(s) = I\{F(y_k) \leq s\} - s$. Our results rely on the following recent theorem, which is a special case of Theorem 2 in Berkes *et al.* [4]. We write $a_n \ll b_n$ if $\limsup_{n \to \infty} |a_n/b_n| < \infty$. A function $F$ is Lipschitz continuous of order $\theta$ if there is some constant $\kappa$ such $|F(y) - F(x)| \leq \kappa |y - x|^\theta$ for all $x, y \in \mathbb{R}$.

**Theorem 8.** *Assume that $\{y_k, k \in \mathbb{Z}\}$ is a strictly stationary sequence which can be represented in the form*

$$y_k = f(\ldots, \varepsilon_{k-1}, \varepsilon_k),$$

*where $f : \mathbb{R}^\mathbb{N} \to \mathbb{R}$ is Borel measurable and $\{\varepsilon_k, k \in \mathbb{Z}\}$ is an i.i.d. sequence. Further, assume that the distribution function $F(x) = P(y_0 \leq x)$ is Lipschitz continuous of order $\theta$. If there are measurable functions $f_m : \mathbb{R}^m \to \mathbb{R}$ ($m \geq 1$) such that*

$$P(|y_k - f_m(\varepsilon_{k-m+1}, \ldots, \varepsilon_k)| > m^{-A}) \ll m^{-B} \qquad \text{for } \min\{A/\theta, B\} > 4,$$

*then the series*

$$\Gamma(s, s') = \sum_{-\infty < k < \infty} \mathbb{E} Y_0(s) Y_k(s') \tag{26}$$

*converges absolutely for every choice of parameters $0 \leq s, s' \leq 1$. Moreover, we can construct $\{y_k, k \in \mathbb{Z}\}$ and a two-parameter Gaussian process $\{K(s,t), (s,t) \in [0,1]^2\}$, with $\mathbb{E} K(s,t) = 0$ and $\mathbb{E} K(s,t) K(s', t') = (t \wedge t') \Gamma(s, s')$, on one probability space such that for some $\alpha > 0$,*

$$\sup_{0 \leq t \leq T} \sup_{0 \leq s \leq 1} |R(s,t) - K(s,t)| = o(T^{1/2} (\log T)^{-\alpha}) \qquad \text{a.s.}$$

A look at our approximation results in Section 3 (Lemmas 6 and 7) shows that Theorem 8 can be almost directly applied for augmented GARCH sequences. The only thing that remains to be checked is the Lipschitz continuity of $F$. This will follow if, for example, we demand that $H(x) = P(\varepsilon_0 \leq x)$ is Lipschitz continuous of order $\theta$ and that

$$\mathbb{E} \sigma_0^{-\theta} < \infty. \tag{27}$$

Indeed, since $\varepsilon_0$ and $\sigma_0$ are independent, we get

$$|P(y_0 \leq x_2) - P(y_0 \leq x_1)| = |P(\varepsilon_0 \leq x_2/\sigma_0) - P(\varepsilon_0 \leq x_1/\sigma_0)|$$
$$\leq \mathbb{E}|H(x_2/\sigma_0) - H(x_1/\sigma_0)|$$
$$\leq L\mathbb{E}\sigma_0^{-\theta}|x_2 - x_1|^\theta.$$

A sufficient condition for (27) is $\sigma_0 \geq \delta > 0$, which is satisfied in all known examples for polynomial GARCH processes. In case of exponential GARCH, the conditions

$$|c(\varepsilon_0)| \leq c < 1 \quad \text{and} \quad \mathbb{E}\exp(\theta/2|g(\varepsilon_0)|) < \infty$$

suffice to ensure (27). To see this, note that

$$\mathbb{E}\sigma_0^{-\theta} \leq \mathbb{E}\exp(\theta/2|\Lambda(\sigma_0^2)|).$$

Hence, we use the same arguments as in the proof of Proposition 1 to show that the last term is finite. The following theorems summarize our observations.

**Theorem 9.** *Let $\Lambda(x) = x^\delta$, $\delta > 0$. Assume that (1), (4) and (7) hold. Assume, further, that the distribution function $H(x) = P(\varepsilon_0 \leq x)$ is Lipschitz continuous of order $\theta > 0$ and that (27) holds. If (20) holds and $\mathbb{E}\log^+|\varepsilon_0|^{(\mu-2)/2} < \infty$ with some $\mu > 10$, then the proposition of Theorem 8 is valid.*

**Theorem 10.** *Let $\Lambda(x) = \log x$. Assume that (1), (4) and (7) hold. Assume, further, that the distribution function $H(x) = P(\varepsilon_0 \leq x)$ is Lipschitz continuous of order $\theta > 0$ and that (27) holds. If for some $\mu > 10$, (21) and $\mathbb{E}\log^+|\varepsilon_0|^{(\mu-2)/2} < \infty$ hold, then the proposition of Theorem 8 is valid.*

An immediate consequence is the following two-parameter version of the empirical central limit theorem which can, for example, be used to detect a change in the structure of the GARCH sequence (cf. Bai [2]).

**Theorem 11.** *Assume that $\{y_k, k \in \mathbb{Z}\}$ is an augmented GARCH sequence satisfying the conditions of either Theorem 9 or Theorem 10. Let $\{K(s,t), (s,t) \in [0,1]^2\}$ be a Gaussian process with $\mathbb{E}K(s,t) = 0$ and $\mathbb{E}K(s,t)K(s',t') = (t \wedge t')\Gamma(s,s')$. Then,*

$$n^{1/2}\left(\frac{1}{n}\sum_{1 \leq k \leq nt}(I\{F(y_k) \leq s\} - s)\right) \overset{\mathcal{D}[0,1]^2}{\longrightarrow} K(s,t) \qquad (n \to \infty).$$

As we pointed out in the introduction, our approach is, in many cases, superior to results following from mixing properties. Carrasco and Chen [11] verified $\beta$-mixing with exponential decay for $y_k$, an approach requiring a continuous positive density of $\varepsilon_0$ on $(-\infty, \infty)$ and $\mathbb{E}\varepsilon_0 = 0$, $\mathbb{E}\varepsilon_0^2 = 1$ and

$$|c(0)| < 1, \qquad \mathbb{E}|c(\varepsilon_0)| < 1, \qquad \mathbb{E}|g(\varepsilon_0)| < \infty. \tag{28}$$

Together with Theorem 2 in Philipp and Pinzur [26], this yields a similar result. Clearly, (28) is more restrictive than (20). Also, our results do not require a positive and continuous density. In the literature, special attention has been paid to the IGARCH(1, 1) process, that is, GARCH(1, 1) with $\mathbb{E}(\beta + \alpha\varepsilon_0^2) = 1$. For an IGARCH process, (28) does not hold, but our new results still apply. The moment conditions (28) for $c(\varepsilon_0)$ and $g(\varepsilon_0)$ are slightly milder than those of Theorem 10 in case of exponential GARCH sequences.

## 5. Proofs

### 5.1. Perturbation error

**Proof of Lemma 1.** We assume, without loss of generality, that $k = 0$. By our assumption (9), the term $\Lambda(\sigma_{km}^2)$ is non-negative and thus, using the special choice of $f$ and $\Lambda$, we get

$$|\eta_0 - \eta_{0m}|^2 = |\varepsilon_0|^{2\nu}|(\Lambda^{-1} \circ \Lambda(\sigma_0^2))^{\nu/2} - (\Lambda^{-1} \circ \Lambda(\sigma_{0m}^2))^{\nu/2}|^2$$

$$\leq |\varepsilon_0|^{2\nu}|(\Lambda^{-1} \circ \Lambda(\sigma_0^2))^{\nu} - (\Lambda^{-1} \circ \Lambda(\sigma_{0m}^2))^{\nu}| \quad (29)$$

$$= |\varepsilon_0|^{2\nu}|\Lambda(\sigma_0^2)^{\nu/\delta} - \Lambda(\sigma_{0m}^2)^{\nu/\delta}|. \quad (30)$$

Let us first consider the case $\nu/\delta \leq 1$. From (7) and Minkowski's inequality (Hardy, Littlewood and Pólya [20]), we infer that

$$|\Lambda(\sigma_0^2)^{\nu/\delta} - \Lambda(\sigma_{0m}^2)^{\nu/\delta}| \leq \left(\sum_{i=m+1}^{\infty} g(\varepsilon_{-i}) \prod_{1 \leq j < i} c(\varepsilon_{-j})\right)^{\nu/\delta}$$

$$\leq \sum_{i=m+1}^{\infty} g(\varepsilon_{-i})^{\nu/\delta} \prod_{1 \leq j < i} c(\varepsilon_{-j})^{\nu/\delta}.$$

Since we assume $\mathbb{E}|\eta_0|^2 < \infty$, it follows from Corollary 1 that

$$\mathbb{E}|\eta_0 - \eta_{0m}|^2 \leq \mathbb{E}|\varepsilon_0|^{2\nu} \sum_{i=m+1}^{\infty} \mathbb{E}g(\varepsilon_{-i})^{\nu/\delta} \prod_{1 \leq j < i} \mathbb{E}c(\varepsilon_{-j})^{\nu/\delta} \leq c_1 \varrho_1^m,$$

with some constant $c_1 > 0$ and $\varrho_1 = \mathbb{E}c(\varepsilon_0)^{\nu/\delta} < 1$.

If $\nu/\delta > 1$, then by the mean value theorem, (30) is bounded by

$$|\varepsilon_0|^{2\nu}\frac{\nu}{\delta}|\Lambda(\sigma_0^2)|^{\nu/\delta-1}|\Lambda(\sigma_0^2) - \Lambda(\sigma_{0m}^2)|.$$

From this, we get, by the Hölder and the Minkowski inequalities, that

$$\mathbb{E}|\eta_0 - \eta_{0m}|^2$$

$$\leq \frac{\nu}{\delta}\mathbb{E}|\varepsilon_0|^{2\nu}(\mathbb{E}|\Lambda(\sigma_0^2)|^{\nu/\delta})^{(\nu-\delta)/\nu}\left(\mathbb{E}\left(\sum_{i=m+1}^{\infty} g(\varepsilon_{-i})\prod_{1\leq j<i}c(\varepsilon_{-j})\right)^{\nu/\delta}\right)^{\delta/\nu}$$

$$\leq c_2 \sum_{i=m+1}^{\infty}\left(\mathbb{E}g(\varepsilon_{-i})^{\nu/\delta}\prod_{1\leq j<i}\mathbb{E}c(\varepsilon_{-j})^{\nu/\delta}\right)^{\delta/\nu} \leq c_3(\varrho_1^{\delta/\nu})^m.$$

□

**Proof of Lemma 2.** Assume, again, that $k=0$. We can formally derive (29), which reads here as

$$|\eta_0 - \eta_{0m}|^2 \leq |\varepsilon_0|^{2\nu}|e^{\nu\log\sigma_0^2} - e^{\nu\log\sigma_{0m}^2}|. \tag{31}$$

By the mean value theorem, we get

$$|e^{\nu\log\sigma_0^2} - e^{\nu\log\sigma_{0m}^2}| \leq \nu(e^{\nu\log\sigma_0^2} + e^{\nu\log\sigma_{0m}^2})|\log\sigma_0^2 - \log\sigma_{0m}^2|. \tag{32}$$

Since $\nu < \mu$, we can find some $\zeta > 0$ such that $\nu(1+\zeta) < \mu$. It follows from Proposition 1 that

$$\mathbb{E}e^{\nu(1+\zeta)\log\sigma_0^2} < \infty$$

and, similarly, we get

$$\mathbb{E}e^{\nu(1+\zeta)\log\sigma_{0m}^2} < \infty.$$

Thus, the Hölder inequality and the Minkowski inequality give

$$\mathbb{E}e^{\nu\log\sigma_0^2}|\log\sigma_0^2 - \log\sigma_{0m}^2|$$
$$\leq (\mathbb{E}e^{\nu(1+\zeta)\log\sigma_0^2})^{1/(1+\zeta)}(\mathbb{E}|\log\sigma_0^2 - \log\sigma_{0m}^2|^{(1+\zeta)/\zeta})^{\zeta/(1+\zeta)}$$
$$\leq (\mathbb{E}e^{\nu(1+\zeta)\log\sigma_0^2})^{1/(1+\zeta)}(\mathbb{E}|g(\varepsilon_0)|^{(1+\zeta)/\zeta})^{\zeta/(1+\zeta)}\sum_{i=m+1}^{\infty}c^{i-1}\leq \text{const}\cdot c^m.$$

The analogous result can be obtained if we replace $\exp(\nu\log\sigma_0^2)$ by $\exp(\nu\log\sigma_{0m}^2)$. Relation (18) now follows from (31) and (32). □

**Proof of Lemma 4.** We can make use of Lemma 3. By the first part of (20) and the Rosenthal inequality (see, e.g., Petrov [25], Theorem 2.9), we infer that $P(T_m > C_3 m) \leq c_1 m^{-\mu/2}$. Hence, it remains to estimate $P(|\Lambda(\sigma_0^2)| > e^{C_4 m})$. We set $\gamma_1 = \mathbb{E}\log|c(\varepsilon_0)|$ and $\gamma_2 = \mathbb{E}\log^+|g(\varepsilon_0)|$. Let $0 \leq \varrho < 1$ such that $\log\varrho - \gamma_1 = a > 0$. If $t$ is chosen sufficiently large, we get, from (7) and the Markov inequality,

$$P(|\Lambda(\sigma_0^2)| > t)$$

$$\leq \sum_{i=1}^{\infty} P\left(|g(\varepsilon_{-i})| \prod_{1\leq j<i} |c(\varepsilon_{-j})| > t(1-\varrho)\varrho^{i-1}\right)$$

$$= \sum_{i=1}^{\infty} P\left(\log^+ |g(\varepsilon_{-i})| + \sum_{1\leq j<i} (\log|c(\varepsilon_{-j})| - \gamma_1) > \log t - c_2 + i(\log\varrho - \gamma_1)\right)$$

$$\leq \sum_{i=1}^{\infty} P\left(\log^+ |g(\varepsilon_{-i})| - \gamma_2 + \sum_{1\leq j<i} (\log|c(\varepsilon_{-j})| - \gamma_1) > \frac{\log t}{2} + ia\right)$$

$$\leq \sum_{i=1}^{\infty} \mathbb{E}\left|\log^+ |g(\varepsilon_{-i})| - \gamma_2 + \sum_{1\leq j<i} (\log|c(\varepsilon_{-j})| - \gamma_1)\right|^\mu \left(\frac{\log t}{2} + ia\right)^{-\mu}.$$

From (20), and again from the Rosenthal inequality, we obtain

$$P(|\Lambda(\sigma_0^2)| > t) \leq c_3 \sum_{i=1}^{\infty} i^{\mu/2} \left(\frac{\log t}{2} + ia\right)^{-\mu}$$

$$\leq c_4 \sum_{i=1}^{\infty} \left(\frac{\log t}{2} + ia\right)^{-\mu/2} \leq c_5 (\log t)^{-(\mu-2)/2}. \qquad \square$$

Since Lemmas 5–7 follow by similar arguments, we will omit their proofs.

### 5.2. A Berry–Esseen bound

We shall prove Theorem 6; the arguments for Theorem 7 are identical. Our proof relies on the following result for $m$-dependent sequences.

**Lemma 8 (Tihomirov [27]).** *Let $X_1, X_2, \ldots$ be a strictly stationary sequence of $m$-dependent random variables with $\mathbb{E}X_1 = 0$ such that*

$$\mathbb{E}|X_1|^3 < \infty.$$

*Let $B_n^2 = \mathrm{Var}(S_n)$ and $\beta^2 = \mathbb{E}X_1^2 + 2\sum_{k=2}^m \mathbb{E}X_1 X_k$. If $\beta^2 > 0$, then there exist absolute constants $C_1$ and $C_2$ such that*

$$\sup_{x\in\mathbb{R}} |P\{S_n \leq B_n x\} - \Phi(x)| \leq C_1 \frac{b_m^2 \mathbb{E}^{1/3}|X_1|^3}{\beta^3 \sqrt{n}} + C_2 \frac{mb_m \mathbb{E}^{1/3}|X_1|^3 \log n}{\beta^2 n},$$

*where $b_m = \max_{1\leq p\leq m+1} \mathbb{E}^{1/3}|\sum_{v=1}^p X_v|^3$.*

We use the notation of Section 4.3 and additionally define $\eta_k$ and $\eta_{km}$ by

$$\eta_k = f(y_k) - \mathbb{E}f(y_k) \quad \text{and} \quad \eta_{km} = f(\varepsilon_k \sigma_{km}) - \mathbb{E}f(\varepsilon_k \sigma_{km}).$$

Then, clearly, for every $m \geq 1$, $\{\eta_{km}, k \in \mathbb{Z}\}$ defines some strictly stationary and $m$-dependent sequence. By a routine argument, we can extend the proof of Lemma 1 to show

$$\|\eta_k - \eta_{km}\|_2 \leq c_1 \varrho^m \qquad \text{for all } m \geq 1 \text{ and } k \in \mathbb{Z}, \tag{33}$$

for some positive constant $c_1$ and $\varrho < 1$. We choose some $n \in \mathbb{N}$ which we assume to be fixed for the moment and set $S'_n = \sum_{k=1}^n \eta_{km}$ and $(B'_n)^2 = \operatorname{Var} S'_n$ with $m = [t \log n]$. The value of $t$ will be specified later. Let $(\delta_k)$ be a sequence of positive reals. A simple estimate gives

$$P\{S_n > xB_n\} \leq P\{S'_n > (x - \delta_n)B_n\} + P\{|S_n - S'_n| > \delta_n B_n\}$$

and

$$P\{S'_n > (x + \delta_n)B_n\} \leq P\{S_n > xB_n\} + P\{|S_n - S'_n| > \delta_n B_n\}.$$

A repeated application of the triangle inequality, together with the latter estimates, shows, after a moment's reflection, that

$$\sup_{x \in \mathbb{R}} |P\{S_n \leq B_n x\} - \Phi(x)| \leq R_n^{(1)} + R_n^{(2)} + R_n^{(3)}, \tag{34}$$

where

$$R_n^{(1)} = \sup_{x \in \mathbb{R}} |P\{S'_n \leq B_n(x + \delta_n)\} - \Phi(x + \delta_n)|,$$

$$R_n^{(2)} = \sup_{x \in \mathbb{R}} |\Phi(x + \delta_n) - \Phi(x)|,$$

$$R_n^{(3)} = P\{|S_n - S'_n| > \delta_n B_n\}.$$

We shall now estimate $R_n^{(i)}$, $i = 1, 2, 3$. Setting $\delta_n = n^{-1/2}$, we get, by the mean value theorem, $R_n^{(2)} \leq (2\pi n)^{-1/2}$. From (25), we conclude that $2B_n^2 > \beta^2 n$ for all $n \geq n_0$. Since we assume that $\beta^2 > 0$, we infer from the Markov and the Minkowski inequalities and (33) that

$$R_n^{(3)} = P\{|S_n - S'_n| > \delta_n B_n\} \leq n\mathbb{E}|S_n - S'_n|^2 / B_n^2$$

$$\leq n \left( \sum_{k=1}^n \|y_k - y_{km}\|_2 \right)^2 \Big/ B_n^2 \leq c_2 n^2 \varrho^m,$$

where $c_2$ can be chosen such that it is not dependent on $n$. (The following constants $c_i$ occurring in the proof are also independent of $n$.) Hence, if $t \geq -5/\log \varrho$, we have $R_n^{(3)} = O(n^{-1/2})$. In order to estimate $R_n^{(1)}$, we will use Lemma 8. This cannot be directly applied; we have to change from $B_n$ to $B'_n$. By Cauchy–Schwarz and Minkowski's inequality again,

we derive

$$|B_n^2 - (B_n')^2| \le \mathbb{E}|S_n - S_n'||S_n + S_n'| \le c_3\|S_n - S_n'\|_2 B_n$$
$$\le c_4 n^{1/2}\sum_{k=1}^n \|y_k - y_{km}\|_2 \le c_5 n^{3/2}\varrho^{m/2} \le c_6 n^{-1},$$

where the last inequality follows from the choice $t \ge -5/\log\varrho$. Again using $2B_n^2 > \beta^2 n$, if $n \ge n_0$, we can reformulate the latter estimate to

$$(1 + c_7 n^{-2})^{-1/2} \le \frac{B_n}{B_n'} \le (1 - c_7 n^{-2})^{-1/2}, \tag{35}$$

provided that $c_7 n^{-2} < 1$. By routine arguments (cf. [25], Lemma 5.2), it follows that

$$\sup_{y\in\mathbb{R}} |\Phi(yp) - \Phi(y)| \le \begin{cases} (p-1)/(2\pi\mathrm{e})^{1/2}, & \text{if } p \ge 1, \\ (p^{-1} - 1)/(2\pi\mathrm{e})^{1/2}, & \text{if } 0 < p < 1. \end{cases} \tag{36}$$

Trivially, we have

$$R_n^{(1)} \le |P\{S_n' \le B_n(x + \delta_n)\} - \Phi(B_n/B_n'(x + \delta_n))|$$
$$+ |\Phi(B_n/B_n'(x + \delta_n)) - \Phi(x + \delta_n)|$$
$$=: R_n^{(11)} + R_n^{(12)}.$$

From (35) and (36), we get $R_n^{(12)} = O(n^{-2})$. Next, observe that $\mathbb{E}|\eta_{1m}|^3 \le c_8\mathbb{E}|f(y_1)|^3$ and that as consequence of (35), we have $\beta_n' \sim \beta$, where

$$\beta_n' = \left(\mathrm{Var}\,\eta_{1m} + 2\sum_{j=2}^{[t\log n]} \mathrm{Cov}(\eta_{1m}, \eta_{jm})\right)^{1/2}.$$

By Lemma 8, it follows that $R_n^{(11)} = O(n^{-1/2}\log^2 n)$. Finally, combining our estimates for $R_n^{(11)}$, $R_n^{(12)}$, $R_n^{(2)}$ and $R_n^{(3)}$ with (34) completes the proof.

## Acknowledgement

Parts of this paper were developed during the author's Ph.D. studies. The author is grateful to his supervisor István Berkes for his valuable comments.